\documentclass{amsart}
\usepackage{amssymb, latexsym} 
\newtheorem{theorem}{Theorem}

\newtheorem{remark}{Remark}
\newtheorem{example}{Example}

\begin{document}
\title{On Products of Closed Subsets in Free Groups}
\author{Rita Gitik}
\address{ Department of Mathematics \\ University of Michigan \\ Ann Arbor, MI, 48109}
\email{ritagtk@umich.edu}
\author{Eliyahu Rips}
\address{ Institute of Mathematics \\ Hebrew University, Jerusalem, 91904, Israel}
\date{\today}

\begin{abstract}
We present examples of closed subsets of a free group such that their product is not closed in the profinite topology. We discuss how to characterize a subset of a free group which is closed in the profinite topology and its product with any finitely generated subgroup of a free group is also closed in the profinite topology.
\end{abstract}


\maketitle


\section{Introduction}

The products of closed subsets of topological groups have been studied for a long time.
The following well-known fact is an important result in the general theory of topological groups.

\begin{theorem} cf. \cite{H-R}, Theorem 4.4.

Let $G$ be a topological group, let $A$ be a compact subset of $G$ and let $B$ be a closed subset of $G$. Then the set $AB$ is closed.
\end{theorem}

It has been known for a long time that if the set $A$ is closed but not compact, then the set $AB$ might fail to be closed.

\begin{example} 

Let $G=\mathbf{R}, A= \mathbf{Z}$, and $B= \pi \cdot \mathbf{Z}$. Both $A$ and $B$ are closed subsets of $G$, but the set $A+B$ is a proper dense subset of $\mathbf{R}$, so it is not closed in $G$. 
\end{example}

Note that in Example 1 both $A$ and $B$ are subgroups of $G$, so the product of closed subgroups is not closed in an arbitrary topological group. The situation is different in the profinite topology, which is defined by
proclaiming all subgroups of finite index of $G$ and their cosets to be basic open sets.

Denote the profinite topology on a group $G$ by $PT(G)$.
An open set in $PT(G)$ is a (possibly infinite) union of cosets of various subgroups of finite index and
the closed sets in $PT(G)$ are the complements of such unions in $G$.

A well-known theorem of M. Hall \cite{Hall} states (in  different language) 
that any finitely generated subgroup of a free group $F$ is closed 
in $PT(F)$. This result has been  generalized by many researchers.
  
The authors proved in \cite{G-R1} and \cite{G-R2} 
that the product of two finitely generated subgroups of a free group $F$
is closed in $PT(F)$. A different proof of that fact was published by G.A. Niblo, cf. \cite{Ni}. 
Note that a finitely generated subgroup of a free group need not be compact, so this result cannot be deduced from Theorem 1.
 
A more general result saying that for any finitely generated subgroups
$H_1, \cdots ,H_n$ of a free group $F$  the product $H_1 H_2 \cdots H_{n-1} H_n$ is closed in  $PT(F)$ was obtained by L. Ribes and P.A. Zalesskii in
\cite{R-Z}, by K. Henckell, S.T. Margolis, J.E. Pin, and J. Rhodes in \cite{H-M-P-R},
and by B. Steinberg in \cite{Ste}.
 
T. Coulbois in \cite{Co} proved that property $RZ_n$ is closed under free products, 
where a group $G$ is said to have property $RZ_n$ if for any $n$ finitely generated subgroups 
$H_1, \cdots  H_n$ of $G$, the product $H_1 H_2 \cdots H_{n-1} H_n$ is closed in $PT(G)$. 

The first author showed in \cite{Gi} that if $H$ and $K$ are quasiconvex subgroups of a negatively curved LERF group $G$ and $H$ is malnormal in $G$ then the double coset $KH$ is closed in $PT(G)$. A. Minasyan generalized that result in \cite{Mi}, showing that a product of any finite number of quasiconvex subgroups of a negatively curved GFERF group $G$ is closed in $PT(G)$.(A negatively curved group $G$ is GFERF if all its quasiconvex subgroups are closed in $PT(G)$.)

The aforementioned results lead to the following question: is it true that for any finitely generated subgroup $H$ of a free
group $F$ and for any subset $S$ of $F$, which is closed in $PT(F)$, the product $SH$ is closed in  $PT(F)$.

As a finitely generated subgroup of $F$ need not be compact in $PT(F)$,  we cannot apply Theorem 1 to resolve this question.
Not surprisingly, the answer to this question is negative, in general, as shown in Example 2, Example 3, and Example 4 in Section 2 of this paper.

So we would ask a more restricted question: can we characterize subsets $S$ of $F$, which are closed in $PT(F)$, such that for any finitely generated subgroup $H$ of $F$ the product $SH$ is closed in $PT(F)$. The following special case might be considered first.

\textbf{Question 1.}

Consider a subset $S$ of a free group $F$ such that $S$ is closed in $PT(F)$. Assume that for any element $c$ of $F$ the set $\langle c \rangle \cdot S$ is closed in $PT(F)$. Is it true that for any elements $a$ and $b$ of $F$ the set $\langle a, b \rangle \cdot F$ is closed in $PT(F)$?

\begin{remark} Note that for any subset $S$ of $F$ which is closed in $PT(F)$  there exists a finitely generated subgroup $H_S$ of $F$ such that $SH_S$ is closed in $PT(F)$. Indeed, take $H_S= \langle 1 \rangle$. 
\end{remark}

\section{The Examples}

Let $F$ be the free group of rank two generated by elements $a$ and $b$.

The following example describes an infinitely generated normal subgroup $N$ of $F$ which is closed in $PT(F)$ such that the double coset 
$\langle a^2 \rangle N$ is not closed in $PT(F)$.

\begin{example}

Consider the Higman's group $G= \langle a, b | b^{-1}a b = a^2 \rangle $. As $G$ is metabelian and linear, it is residually finite, however, the cyclic subgroup $\langle a^2 \rangle$ is not separable from the element $a \in G$ in $PT(G)$. Indeed, $a=b^n a^{2^n} b^{-n}$ (which can be shown by induction) and a finite index subgroup of $G$ should contain an element $b^n \in G$ for some positive $n$.  
 
Let $N$ be the kernel of the quotient map from $F= \langle a,b \rangle$ to $G$. As $G$ is residually finite and $N$ is the preimage of the trivial subgroup of $G$ in $F$, it follows that $N$ is closed in $PT(F)$. The subgroup $\langle a^2 \rangle$ of $F$ is closed in
$PT(F)$ because it is finitely generated. However, the double coset $\langle a^2 \rangle N$ is not closed in $PT(F)$ because it is not separable from the element $a \in F$.

Indeed, let $F_0$ be any finite index subgroup of $F$ such that $\langle a^2 \rangle N  \subset F_0$. Let $G_0 < G$ be the projection of $F_0$ to $G$. As 
$a^2 \in G_0$ and $G_0$ has a finite index in $G$, it follows that $a \in G_0$, hence there exists an element $n_0 \in N$ such that $an_0 \in F_0$. However, $a^2 n \in N$ for all $n \in N$, hence $(n_0^{-1}a^{-1})(a^2  n) \in F_0$ for all $n \in N$. As $N$ is normal in $N$, there exists 
$n_1 \in N$ such that $n_0^{-1}a=an_1$. Then $an_1n \in F_0$ for all $n \in N$, hence $aN \subset F_0$. It follows that $a \in F_0$. 
\end{example}

The following example describes a set $S$ which is closed in $PT(F)$ such that its product with 
a free factor of $F$ is not closed in $PT(F)$.

\begin{example} cf. \cite{G-R3}.

Let $F = \langle a,b \rangle$ be the free group of rank two.
Consider an infinite sequence $A = \{ a, a^{2!}, a^{3!}, \cdots, a^{k!}, \cdots \} \subset F$.
Note that $A$ converges to $1_F$. Indeed, let $M$ be a normal subgroup of finite index $m$ in $F$.
If $k \ge m$  then $a^{k!}$ is contained in $M$. 
Hence any open neighborhood of $1_F$ in $F$ contains all, but finitely many elements of $A$,
therefore $A$ converges to $1_F$. Note that $1_F \notin A$, so $A$ is not closed in $PT(F)$. 

Let $m_k, k \ge 1$ be integers such that 
$m_k \rightarrow m_0 \in \hat{Z} \setminus Z$ in $PT(\hat{Z})$,
where $\hat{Z}$ is the completion of $Z$ in $PT(Z)$.
Then $a^k b^{m_k} \rightarrow a^0 b^{m_0} \in \hat{F} \setminus F$,
where $\hat{F}$ is the completion of $F$ in $PT(F)$.
Hence the sequence $a^k b^{m_k}$ has no other limit points.
In particular, it has no limit points in $F$. 
Therefore for every $w \in F$ with 
$w \neq a^k b^{m_k}$ for all $k \ge 1$, there exists an open neighborhood $U$ of $w$
such that $a^k b^{m_k} \notin U$, for all $k \ge 1$.
It follows that the set 
$S= \{ a b ^{m_1}, a^{2!} b ^{m_2}, \cdots , a^{k!} b ^{m_k}, \cdots \}$ is closed in $PT(F)$.

Note that $1_F \notin S \langle b \rangle$, however $ A \subseteq S \langle b \rangle$, so
$1_F \in \bar A \subseteq \overline{S \langle b \rangle}$.
We conclude that $S \langle b \rangle$ is not closed in $PT(F)$.
\end{example}
   
Example 3 motivates the following question.

\textbf{Question 2.}

Is it possible to impose some restrictions on a set $S$ which is closed in $PT(F)$,
such that the product of $S$ with a free factor of $F$ would be closed in $PT(F)$?

\bigskip

The following example demonstrates that such restrictions on $S$ should be severe.

\begin{example} cf. \cite{G-R3}.

Let $F$ be a finitely generated free group on free generators $K \cup L$ such that
$F = \langle K \rangle * \langle L \rangle$. There exists a discrete subset $S$ of $F$ which is closed in $PT(F)$ such that $S \langle K \rangle$ 
is not closed in $PT(F)$ and the last syllable of all elements of $S$ is in $ \langle L \rangle$.
\end{example}

\section{Acknowledgment}

The first author would like to thank the Institute of Mathematics
of the Hebrew University for generous support.

\end{document}